\documentclass[11pt]{article}
\usepackage{epic,latexsym,amssymb}
\usepackage{color}
\usepackage{tikz}
\usepackage{amsmath}

\newtheorem{theorem}{Theorem}

\newtheorem{observation}{Observation}
\newtheorem{claim}{Claim}
\newtheorem{subclaim}{Claim}[claim]
\newtheorem{corollary}{Corollary}

\newcommand{\QED}{$\Box$}
\newcommand{\smallqed}{{\tiny ($\Box$)}}

\newcommand{\cH}{{\cal H}}

\newcommand{\pn}{{\rm pn}}

\newcommand{\diam}{\mathop{\rm diam}}

\newcommand{\proof}{\noindent\textbf{Proof. }}
\newcommand{\2}{ \vspace{0.2cm} }
\newcommand{\1}{ \vspace{0.1cm} }

\let\oldenumerate\enumerate
\renewcommand{\enumerate}{
  \oldenumerate
  \setlength{\itemsep}{0pt}
  \setlength{\parskip}{0pt}
  \setlength{\parsep}{0pt}
}

\begin{document}

\title{Common Domination Perfect Graphs}

\author{$^1$Magda Dettlaff,
$^2$Michael A. Henning\thanks{Research supported in part by the University of Johannesburg} \,
and $^3$Jerzy Topp \\
\\
$^1$Faculty of Mathematics, Physics and Informatics\\
University of Gda\'nsk \\
80-952 Gda\'nsk, Poland \\
\small \tt Email:  magda.dettlaff1@pg.edu.pl \\
\\
$^2$Department of Mathematics and Applied Mathematics\\
University of Johannesburg \\
Auckland Park, 2006, South Africa \\
\small \tt Email: mahenning@uj.ac.za \\
\\
$^3$Institute of Applied Informatics\\
University of Applied Sciences\\
82-300 Elbl\c{a}g, Poland \\
\small \tt Email: j.topp@ans-elblag.pl
}

\date{}
\maketitle

\begin{abstract}
A dominating set in a graph $G$ is a set $S$ of vertices such that every vertex that does not belong to $S$ is adjacent to a vertex in $S$. The domination number $\gamma(G)$ of $G$ is the minimum cardinality of a dominating set of $G$. The common independence number $\alpha_c(G)$ of $G$ is the greatest integer $r$ such that every vertex of $G$ belongs to some independent set of cardinality at least~$r$. The common independence number is squeezed between the independent domination number $i(G)$ and the independence number $\alpha(G)$ of $G$, that is, $\gamma(G) \le i(G) \le \alpha_c(G) \le \alpha(G)$. A graph $G$ is domination perfect if $\gamma(H) = i(H)$ for every induced subgraph $H$ of $G$. We define a graph $G$ as common domination perfect if $\gamma(H) = \alpha_c(H)$ for every induced subgraph $H$ of $G$. We provide a characterization of common domination perfect graphs in terms of ten forbidden induced subgraphs.
\end{abstract}


\section{Introduction}
\label{Intro}

A \emph{dominating set} of a graph $G$ is a set $S$ of vertices of $G$ such that every vertex not in $S$ has a~neighbor in $S$, where two vertices are neighbors in $G$ if they are adjacent. If, in addition, $S$ is an independent set, then $S$ is an \emph{independent dominating set}, abbreviated ID-set, of $G$. The domination number $\gamma(G)$ of $G$ is the minimum cardinality of a dominating set in $G$, while the \emph{independent domination number} $i(G)$ of $G$ is the minimum cardinality of an ID-set in $G$. The \emph{independence number} $\alpha(G)$ of $G$ is the maximum cardinality of an independent set in $G$. A $\gamma$-\emph{set} of $G$ is a dominating set of $G$ of minimum cardinality~$\gamma(G)$, and an $i$-\emph{set} of $G$ is an ID-set of $G$ of minimum cardinality~$i(G)$.  Further, an $\alpha$-\emph{set} of $G$ is an independent set of $G$ of maximum cardinality~$\alpha(G)$.
We refer the reader to the survey \cite{GoddardHenning} of results on independent domination in graphs, and to~\cite{HaHeHe-20,HaHeHe-21,HaHeHe-22} for recent books on domination in graphs.

Dettlaff et al.~\cite{DeLeTo-21} introduced the concept of the \emph{common independence number} of a graph $G$, denoted by $\alpha_c(G)$, as the greatest integer $r$ such that every vertex of $G$ belongs to some independent set in $G$ of cardinality at least~$r$. Thus, the common independence number of $G$ refers to numbers of mutually independent vertices of $G$ and it emphasizes the notion of the individual independence of a vertex of $G$ from other vertices of $G$. The common independence number $\alpha_c(G)$ of $G$ is the limit of symmetry in $G$ with respect to the fact that each vertex of $G$ belongs to an independent set of cardinality $\alpha_c(G)$ in $G$, and there are vertices in $G$ that do not belong to any larger independent set in $G$.

As observed in~\cite{DeLeTo-21}, the common independence number is squeezed between the independent domination number and the independence number of $G$, yielding the following ``domination inequality chain''.

\begin{observation}{\rm (\cite{DeLeTo-21})}
\label{ob:relate}
For every graph $G$ we have $\gamma(G) \le i(G) \le \alpha_c(G) \le \alpha(G)$.
\end{observation}

Motivated by the concept of perfect graphs in the chromatic sense, Sumner and Moore~\cite{SuMo-79} in 1979 defined a graph $G$ to be \emph{domination perfect} if $\gamma(H) = i(H)$ for every induced subgraph $H$ of $G$. Domination perfect graphs are now very well studied in the literature. One of the earliest results on domination perfect graphs is due to Mitchell and Hedetniemi~\cite{MiHe-77} in 1977 who proved that line graphs of trees are domination perfect. In 1978 Allan and Laskar~\cite{AlLa-78} generalized their result and showed that every claw-free graph is domination perfect (see also \cite{Gupta-Singh-Arumugam} for recent generalizations of that result). Some of the pioneering work on domination perfect graphs was done by, among others, Topp and Volkmann~\cite{ToVo-91}, Zverovich and Zverovich~\cite{ZvZv-91,ZvZv-95}, Fulman~\cite{Fu-93}, and Rautenbach and Zverovich~\cite{RaZv-01}.
We refer the reader to the good (but already old) survey on domination perfect (and domination critical) graphs by Sumner~\cite{Su-90}.

The most significant contribution to the study of domination perfect graphs is the 1995 result due to Zverovich and Zverovich~\cite{ZvZv-95} that provided a characterization of domination perfect graphs in terms of seventeen forbidden induced subgraphs. Their result built on earlier work, including important results in~\cite{Fu-93,Su-90,ToVo-91,ZvZv-91}. We remark that in~\cite{CaPl-17} counterexamples were presented to this 1995 characterization of domination perfect graphs due to Zverovich and Zverovich, but these ``counterexamples''  are themselves incorrect, as is the ``characterization'' given in~\cite{CaPl-17}.

Generalizing the concept of domination perfect graphs, given any two graph parameters $\lambda$ and $\mu$ for which $\lambda(G) \le \mu(G)$, a graph $G$ is defined in~\cite{DoGoHeMy-06,Fischermann-Volkmann-Zverovich,GoHe-04,HaHeHe-22} to be \emph{$(\lambda, \mu)$-perfect} if $\lambda(H) = \mu(H)$ for every induced subgraph $H$ of $G$. In particular, a domination perfect graph is a $(\gamma,i)$-perfect graph. In this paper we study $(\gamma,\alpha_c)$-perfect graphs, that is, we study graphs $G$ satisfying $\gamma(H) = \alpha_c(H)$ for every induced subgraph $H$ of $G$. We call such graphs \emph{common domination perfect graphs} since here we study perfect graphs with respect to the domination number and the common independence number. We note that by Observation~\ref{ob:relate} every $(\gamma,\alpha_c)$-perfect graph is a $(\gamma,i)$-perfect graph.

\subsection{Notation}

For notation and graph theory terminology, we in general follow~\cite{HaHeHe-22}. Specifically, let $G$ be a graph with vertex set $V(G)$ and edge set $E(G)$, and of order~$n(G) = |V(G)|$ and size $m(G) = |E(G)|$. Two vertices in $G$ are \emph{neighbors} if they are adjacent. The \emph{open neighborhood} $N_G(v)$ of a vertex $v$ in $G$ is the set of neighbors of $v$, while the \emph{closed neighborhood} of $v$ is the set $N_G[v] = \{v\} \cup N_G(v)$. For a set $S \subseteq V(G)$, its \emph{open neighborhood} is the set $N_G(S) = \bigcup_{v \in S} N_G(v)$, and its \emph{closed neighborhood} is the set $N_G[S] = N_G(S) \cup S$. For a set of vertices $S \subseteq V(G)$, the subgraph induced by $S$ is denoted by $G[S]$. For a set $S \subseteq V(G)$ and a vertex $v \in S$, the \emph{$S$-private neighborhood} $\pn[v,S]$ of $v$ is the set of vertices that are in the closed neighborhood of $v$ but not in the closed neighborhood of the set $S \setminus \{v\}$, that is, $\pn[v,S] = \{w \in V \colon N_G[w] \cap S = \{v\}\}$. If $\pn[v,S] \ne \emptyset$, then a vertex in $\pn[v,S]$ is called an \emph{$S$-private neighbor of $v$}.

A graph $G$ is said to be \emph{chordal} if every cycle of $G$ of length four or more contains a~\emph{chord}, i.e., an edge joining two non-consecutive vertices of the cycle. A \emph{block} in a~graph $G$ is a~maximal connected subgraph having the property that it contains no cut-vertex of its own. An \emph{end-block} of $G$ is a block containing exactly one cut-vertex, while an \emph{inner-block} of $G$ is a block containing at least two cut-vertices. A graph is a \emph{block graph} if every block of $G$ is a~complete graph. It is immediately obvious that every block graph is a chordal graph.

The \emph{line graph} of a graph $H$, denoted $L(H)$, is a graph in which vertices are the edges of $H$, with two vertices of $L(H)$ adjacent whenever the corresponding edges of $H$ are. A graph $G$ is a \emph{line graph} if $G=L(H)$ for some graph $H$. The \emph{middle graph} $M(H)$ of a graph $H$ is the line graph of the corona $H \circ K_1$, that is, $M(H) = L(H \circ K_1)$. The \emph{total graph} $T(H)$ of a graph $H$ is the graph with vertex set $V(H) \cup E(H)$ in which two vertices $x$ and $y$ are adjacent if and only if either $x$ and $y$ are adjacent vertices of $H$, or $x$ and $y$ are adjacent edges of $H$, or $x$ is a vertex of $H$ and $y$ is an edge of $H$ incident with $x$.

If $G$ and $H$ are graphs, then the graph $G$ is $H$-\emph{free} if $G$ does not contain the graph $H$ as an induced subgraph. In particular, a \emph{claw}-\emph{free graph} is a graph that is $K_{1,3}$-free. Further, if $\cH$ is a set of graphs, then $G$ is $\cH$-\emph{free} if $G$ does not contain the graph in $\cH$ as an induced subgraph.  We denote a \emph{path}, a \emph{cycle}, and a \emph{complete} graph on $n$ vertices by $P_n$, $C_n$, and $K_n$, respectively. For an integer $k \ge 1$, we let $[k] = \{1,2,\ldots,k\}$.

\section{Main result}

\begin{figure}[htb]
\begin{center}
\begin{tikzpicture}[scale=.8,style=thick,x=0.8cm,y=0.8cm]
\def\vr{2.5pt} 
\path (0,0) coordinate (u1);
\path (0,1.5) coordinate (u2);
\path (0.5,0.75) coordinate (u3);
\path (1.5,0.75) coordinate (u4);
\path (2,0) coordinate (u5);
\path (2,1.5) coordinate (u6);
%
\draw (u2)--(u1)--(u3);
\draw (u4)--(u5)--(u6);
\draw (u1)--(u5);

\draw (u1) [fill=white] circle (\vr);
\draw (u2) [fill=white] circle (\vr);
\draw (u3) [fill=white] circle (\vr);
\draw (u4) [fill=white] circle (\vr);
\draw (u5) [fill=white] circle (\vr);
\draw (u6) [fill=white] circle (\vr);
\draw (1,-0.8) node {{\small (a) $H_1$}};
\path (4,0) coordinate (v1);
\path (4,1.5) coordinate (v2);
\path (4.5,0.75) coordinate (v3);
\path (5.5,0.75) coordinate (v4);
\path (6,0) coordinate (v5);
\path (6,1.5) coordinate (v6);
%
\draw (v2)--(v1)--(v3);
\draw (v4)--(v5)--(v6);
\draw (v1)--(v5);
\draw (v3)--(v4);
\draw (v1) [fill=white] circle (\vr);
\draw (v2) [fill=white] circle (\vr);
\draw (v3) [fill=white] circle (\vr);
\draw (v4) [fill=white] circle (\vr);
\draw (v5) [fill=white] circle (\vr);
\draw (v6) [fill=white] circle (\vr);
\draw (5,-0.8) node {{\small (b) $H_2$}};
\path (8,0) coordinate (w1);
\path (8,1.5) coordinate (w2);
\path (8.5,0.75) coordinate (w3);
\path (9.5,0.75) coordinate (w4);
\path (10,0) coordinate (w5);
\path (10,1.5) coordinate (w6);
%
\draw (w2)--(w1)--(w3);
\draw (w4)--(w5)--(w6);
\draw (w1)--(w5);
\draw (w3)--(w4);
\draw (w3)--(w6);
\draw (w1) [fill=white] circle (\vr);
\draw (w2) [fill=white] circle (\vr);
\draw (w3) [fill=white] circle (\vr);
\draw (w4) [fill=white] circle (\vr);
\draw (w5) [fill=white] circle (\vr);
\draw (w6) [fill=white] circle (\vr);
\draw (9,-0.8) node {{\small (c) $H_3$}};
\path (12,0) coordinate (x1);
\path (12,1.5) coordinate (x2);
\path (12.5,0.75) coordinate (x3);
\path (13.5,0.75) coordinate (x4);
\path (14,0) coordinate (x5);
\path (14,1.5) coordinate (x6);
%
\draw (x2)--(x1)--(x3);
\draw (x4)--(x5)--(x6);
\draw (x1)--(x5);
\draw (x3)--(x4);
\draw (x3)--(x6);
\draw (x2)--(x4);
\draw (x2)--(x6);
\draw (x1) [fill=white] circle (\vr);
\draw (x2) [fill=white] circle (\vr);
\draw (x3) [fill=white] circle (\vr);
\draw (x4) [fill=white] circle (\vr);
\draw (x5) [fill=white] circle (\vr);
\draw (x6) [fill=white] circle (\vr);
\draw (13,-0.8) node {{\small (d) $H_4$}};
\path (16,0) coordinate (y1);
\path (16,1.5) coordinate (y2);
\path (16.5,0.75) coordinate (y3);
\path (17.5,0.75) coordinate (y4);
\path (18,0) coordinate (y5);
\path (18,1.5) coordinate (y6);
%
\draw (y2)--(y1)--(y3);
\draw (y4)--(y5)--(y6);
\draw (y1)--(y5);
\draw (y3)--(y4);
\draw (y2)--(y6);
\draw (y1) [fill=white] circle (\vr);
\draw (y2) [fill=white] circle (\vr);
\draw (y3) [fill=white] circle (\vr);
\draw (y4) [fill=white] circle (\vr);
\draw (y5) [fill=white] circle (\vr);
\draw (y6) [fill=white] circle (\vr);
\draw (17,-0.8) node {{\small (e) $H_5$}};
\path (20,0) coordinate (z1);
\path (20,1.5) coordinate (z2);
\path (20.5,0.75) coordinate (z3);
\path (21.5,0.75) coordinate (z4);
\path (22,0) coordinate (z5);
\path (22,1.5) coordinate (z6);
%
\draw (z2)--(z1)--(z3);
\draw (z4)--(z5)--(z6);
\draw (z1)--(z5);
\draw (z3)--(z4);
\draw (z3)--(z6);
\draw (z2)--(z6);
\draw (z1) [fill=white] circle (\vr);
\draw (z2) [fill=white] circle (\vr);
\draw (z3) [fill=white] circle (\vr);
\draw (z4) [fill=white] circle (\vr);
\draw (z5) [fill=white] circle (\vr);
\draw (z6) [fill=white] circle (\vr);
\draw (21,-0.8) node {{\small (f) $H_6$}};
\end{tikzpicture}

\vskip 0.25 cm

\begin{tikzpicture}[scale=.8,style=thick,x=0.8cm,y=0.8cm]
\def\vr{2.5pt} 
\path (0,0) coordinate (a1);
\path (0,1) coordinate (a2);
\path (1,0) coordinate (a3);
\path (1,1) coordinate (a4);
\path (2,0) coordinate (a5);
\path (2,1) coordinate (a6);
%
\draw (a1)--(a3)--(a5);
\draw (a2)--(a4)--(a6);
\draw (a1) [fill=white] circle (\vr);
\draw (a2) [fill=white] circle (\vr);
\draw (a3) [fill=white] circle (\vr);
\draw (a4) [fill=white] circle (\vr);
\draw (a5) [fill=white] circle (\vr);
\draw (a6) [fill=white] circle (\vr);
\draw (1,-0.8) node {{\small (g) $H_7$}};
\path (4.5,0) coordinate (b1);
\path (4.5,1) coordinate (b2);
\path (5.5,0) coordinate (b3);
\path (5.5,1) coordinate (b4);
\path (6.5,0) coordinate (b5);
\path (6.5,1) coordinate (b6);
%
\draw (b1)--(b3)--(b5);
\draw (b2)--(b4)--(b6);
\draw (b1)--(b2);
\draw (b1) [fill=white] circle (\vr);
\draw (b2) [fill=white] circle (\vr);
\draw (b3) [fill=white] circle (\vr);
\draw (b4) [fill=white] circle (\vr);
\draw (b5) [fill=white] circle (\vr);
\draw (b6) [fill=white] circle (\vr);
\draw (5.5,-0.8) node {{\small (h) $H_8$}};
\path (9,0) coordinate (c1);
\path (9,1) coordinate (c2);
\path (10,0) coordinate (c3);
\path (10,1) coordinate (c4);
\path (11,0) coordinate (c5);
\path (11,1) coordinate (c6);
%
\draw (c1)--(c3)--(c5);
\draw (c2)--(c4)--(c6);
\draw (c1)--(c2);
\draw (c5)--(c6);
\draw (c1) [fill=white] circle (\vr);
\draw (c2) [fill=white] circle (\vr);
\draw (c3) [fill=white] circle (\vr);
\draw (c4) [fill=white] circle (\vr);
\draw (c5) [fill=white] circle (\vr);
\draw (c6) [fill=white] circle (\vr);
\draw (10,-0.8) node {{\small (i) $H_9$}};
\path (13.5,0) coordinate (d1);
\path (13.5,1) coordinate (d2);
\path (14.5,0) coordinate (d3);
\path (14.5,1) coordinate (d4);
\path (15.5,0) coordinate (d5);
\path (15.5,1) coordinate (d6);
%
\draw (d3)--(d5);
\draw (d2)--(d4)--(d6);
\draw (d1)--(d2);
\draw (d5)--(d6);
\draw (d3)--(d4);
\draw (d1) [fill=white] circle (\vr);
\draw (d2) [fill=white] circle (\vr);
\draw (d3) [fill=white] circle (\vr);
\draw (d4) [fill=white] circle (\vr);
\draw (d5) [fill=white] circle (\vr);
\draw (d6) [fill=white] circle (\vr);
\draw (14.5,-0.8) node {{\small (j) $H_{10}$}};
\end{tikzpicture}
\end{center}
\vskip -0.5 cm
\caption{Minimal common domination imperfect graphs $H_1, \ldots, H_{10}$}
\label{fig:forbidden-graphs}
\end{figure}
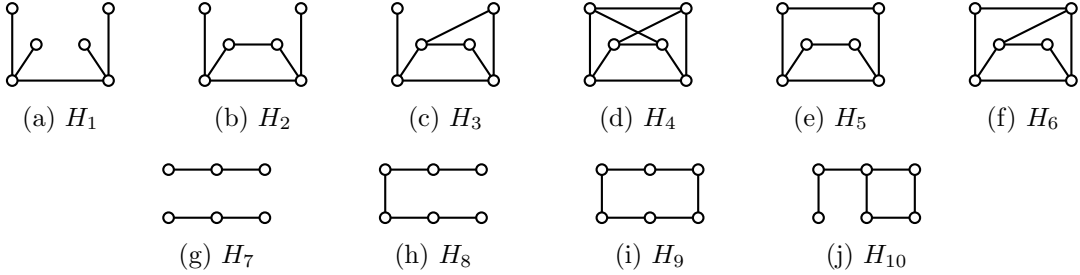

In this paper we provide a characterization of common domination perfect graphs in terms of ten forbidden induced subgraphs. We define a graph $G$ to be \emph{minimal common domination imperfect} if $\gamma(G) < \alpha_c(G)$, yet $\gamma(H) = \alpha_c(H)$ for all proper induced subgraphs $H$ of $G$. Let $\cH = \{H_1, H_2, \ldots, H_{10}\}$, where $H_1, \ldots, H_{10}$ are the ten graphs shown in Figure~\ref{fig:forbidden-graphs}.

We shall need the following properties of the graphs in the family $\cH$.

\begin{observation}
\label{ob:family}
The following properties hold for graphs in the family $\cH$. \\[-24pt]
\begin{enumerate}
\item[{\rm (a)}] $\gamma(H_k) = 2 < 3 = i(H_k) = \alpha_c(H_k)$ if $k \in [4]$.
\item[{\rm (b)}] $\gamma(H_k) = i(H_k) = 2 < 3 = \alpha_c(H_k)$ for $k\in [10] \setminus [4]$.
\item[{\rm (c)}] Every graph in $\cH$ is a minimal common domination imperfect graph.
\end{enumerate}
\end{observation}

In our main theorem we present a forbidden subgraph characterization of $(\gamma,\alpha_c)$-perfect graphs.

\begin{theorem}
\label{gamma-alpha_c-perfect-graphs}
A graph $G$ is a $(\gamma,\alpha_c)$-perfect graph if and only if it does not contain any of the ten graphs $H_1, \ldots, H_{10}$ shown in Figure~\ref{fig:forbidden-graphs} as an induced subgraph.
\end{theorem}
\proof By Observation~\ref{ob:family}, if $G$ is a $(\gamma,\alpha_c)$-perfect graph, then it does not contain any of the ten graphs $H_1, \ldots, H_{10}$ shown in Figure~\ref{fig:forbidden-graphs} as an induced subgraph.

Now assume that $G$ is a graph that does not contain any of the ten graphs $H_1, \ldots, H_{10}$ shown in Figure~\ref{fig:forbidden-graphs} as an induced subgraph. We claim that $G$ is $(\gamma,\alpha_c)$-perfect, that is, $\gamma(H)=\alpha_c(H)$ for every induced subgraph $H$ of $G$. Suppose, to the contrary, that $G$ is a counterexample of minimum order. Thus, none of the graphs $H_1, \ldots, H_{10}$ is an induced subgraph of $G$, yet $\gamma(H) < \alpha_c(H)$ for some induced subgraph $H$ of $G$.

By the minimality of $G$, the graph $G$ is a minimal common domination imperfect graph, and so $\gamma(G) < \alpha_c(G)$ and $\gamma(H) = \alpha_c(H)$ for all proper induced subgraphs $H$ of $G$. If $\gamma(G) = 1$, then $\alpha_c(G) = 1$, contradicting the fact that $\gamma(G) < \alpha_c(G)$. Hence, $\gamma(G) \ge 2$. By Observation~\ref{ob:relate}, $\gamma(G) \le i(G) \le \alpha_c(G)$, implying that either $\gamma(G) < i(G) \le \alpha_c(G)$ or $\gamma(G) = i(G) < \alpha_c(G)$. We proceed further with the following claim.

\begin{claim}
\label{c:claim1}
$\gamma(G) = i(G) < \alpha_c(G)$.
\end{claim}
\proof Suppose, to the contrary, that $\gamma(G) < i(G) \le \alpha_c(G)$. In this case, no $\gamma$-set of $G$ is independent. Among all $\gamma$-sets of $G$, let $D$ be chosen so that the number of edges in $G[D]$ is as small as possible. Since $\gamma(G) < i(G)$, the subgraph $G[D]$ has at least one edge. Let $v_1$ and $v_2$ be two adjacent vertices in $D$, and let $I_{v_1}$ and $I_{v_2}$ be the sets of $D$-private neighbors of $v_1$ and $v_2$, that is, $I_{v_i} = \pn[v_i,D] = \{x\in V(G) \setminus D \colon N_G(x) \cap D = \{v_i\}\}$ for $i \in [2]$. Since $v_1$ and $v_2$ are adjacent and $D$ is a $\gamma$-set of $G$, the sets $I_{v_1}$ and $I_{v_2}$ are nonempty and, obviously, disjoint. Further, $v_i \notin I_{v_i}$ for $i \in [2]$.

If a vertex $y_i \in I_{v_i}$ dominates the set $I_{v_i}$ for some $i \in [2]$, then $I_{v_i} \subseteq N_G[y_i]$ and the set $D' = (D\setminus \{v_i\})\cup \{y_i\}$ would be a $\gamma$-set of $G$ with $m(G[D']) < m(G[D])$, contradicting the choice of the set $D$. Hence, no vertex in $I_{v_i}$ dominates the set $I_{v_i}$. Let $x_i$ and $y_i$ be nonadjacent vertices in $I_{v_i}$ for $i \in [2]$. Now, depending on the possible edges joining the vertices $x_1$ and $y_1$ with the vertices $x_2$ and $y_2$, the induced subgraph $G[\{v_1,x_1,y_1,v_2,x_2,y_2\}]$ of $G$ is one of the forbidden graphs $H_1, \ldots, H_6$ shown in Figure~\ref{fig:forbidden-graphs}, contradicting our assumption that $G$ is $\cH$-free.~\smallqed

\medskip
By Claim~\ref{c:claim1}, $\gamma(G) = i(G) < \alpha_c(G)$. Since $H_8$ is a forbidden induced subgraph of $G$, the graph $G$ is $P_6$-free, implying that the diameter of $G$ is at most~$4$, that is, $\diam(G) \le 4$ and so $d_G(x,y)\le 4$ for every two vertices $x$ and $y$ of $G$. Among all $\gamma$-sets of $G$, let $D$ be chosen to contain the minimum number of vertices whose closed neighborhood induces a complete subgraph of $G$.

\begin{claim}
\label{c:claim2}
$d_G(x,y) \le 3$ for every two vertices $x, y \in D$.
\end{claim}
\proof Suppose, to the contrary, that the $\gamma$-set $D$ of $G$ contains two vertices at distance~$4$ in $G$. Let $v_0$ and $v_4$ be two vertices that belong to $D$ at distance~$4$ in $G$. 
Let $V_i$ be the set of vertices at distance~$i$ from~$v_0$ in $G$ and at distance~$4-i$ from~$v_4$ in $G$ for $i \in [3]$. Thus if $v \in V_i$, then $d_G(v,v_0) = i$ and $d_G(v,v_4) = 4-i$ for $i \in [3]$. The sets $V_1$, $V_2$, $V_3$ are nonempty and disjoint subsets of $V(G)$, and every vertex in $V_2$ is dominated by a vertex belonging to $D\setminus\{v_0, v_4\}$. Let $v_0v_1v_2v_3v_4$ be a shortest $(v_0,v_4)$-path in $G$, and so $v_i \in V_i$ for $i \in [3]$.

\begin{subclaim}
\label{c:claim2.1}
The closed neighborhood of $v_0$ and $v_4$ both induce complete subgraphs of $G$.
\end{subclaim}
\proof We show firstly that both $G[V_1]$ and $G[V_3]$ are complete subgraph. Suppose, to the contrary, that $G[V_1]$ or $G[V_3]$ is not a complete subgraph. By symmetry, and renaming vertices if necessary, we may assume that $G[V_1]$ is not a complete subgraph. Let $x_1$ and $y_1$ be vertices in $V_1$ that are not adjacent. We now consider the induced subgraph $H = G[\{v_0,x_1,y_1,v_2,v_3,v_4\}]$. If neither $x_1$ nor $y_1$ is adjacent to $v_2$, then $H = H_7$. If exactly one of $x_1$ and $y_1$ is adjacent to $v_2$, then $H = H_8$. If both $x_1$ and $y_1$ are adjacent to $v_2$, then $H = H_{10}$. In all three cases, $H \in \cH$, a contradiction. Hence, both $G[V_1]$ and $G[V_3]$ are complete subgraph.

We show next that every neighbor of $v_0$ not in $V_1$ is adjacent to every vertex in $V_1$, and every neighbor of $v_4$ not in $V_3$ is adjacent to every vertex in $V_3$. Suppose this is not the case. By symmetry, and renaming vertices if necessary, we may assume that $v_0$ has a~neighbor $x_1$ that does not belong to $V_1$ and is not adjacent to some vertex $y_1 \in V_1$. Since $x_1 \in N_G(v_0)$ and $x_1 \notin V_1$, we note that the vertex $x_1$ has no neighbor in $V_2$. As before we let $H = G[\{v_0,x_1,y_1,v_2,v_3,v_4\}]$ and deduce that $H = H_7$ or $H = H_8$, a contradiction. Therefore, every vertex belonging to $N_G(v_0) \setminus V_1$ is adjacent to every vertex in $V_1$. Analogously, every vertex belonging to $N_G(v_4) \setminus V_3$ is adjacent to every vertex in $V_3$.

Suppose that $v_0$ has two neighbors $x_1$ and $y_1$ that do not belong to $V_1$ and are not adjacent. In this case, the induced subgraph $G[\{v_0,x_1,y_1,v_2,v_3,v_4\}]$ is isomorphic to $H_7$, a~contradiction. Hence, the subgraph $G[N_G[v_0] \setminus V_1]$ is a complete graph. Analogously, the subgraph $G[N_G[v_4] \setminus V_3]$ is a complete graph. From our earlier properties, we therefore infer that $G[N_G[v_0]]$ and $G[N_G[v_4]]$ are complete subgraphs of $G$.~\smallqed

\medskip
By Claim~\ref{c:claim2.1}, $G[N_G[v_0]]$ and $G[N_G[v_4]]$ are complete subgraphs of $G$, that is, the closed neighborhood of $v_0$ induces a complete subgraph of $G$, as does the closed neighborhood of $v_4$. Replacing the vertex $v_0$ in the set $D$ by the vertex $v_1$ produces a new $\gamma$-set $D' = (D \setminus \{v_0\}) \cup \{v_1\}$ of $G$ noting that $N_G[v_0] \subset N_G[v_1]$ and $|D'| = |D| = \gamma(G)$. However, since the two neighbors $v_0$ and $v_2$ of $v_1$ are not adjacent, the closed neighborhood of $v_1$ does not induce a complete subgraph of $G$. Therefore the $\gamma$-set $D'$ contains fewer vertices whose closed neighborhood induces a complete subgraph of $G$ than does the set $D$, contradicting our choice of the set $D$.~\smallqed

\medskip
By Claim~\ref{c:claim2}, $d_G(x,y) \le 3$ for every two vertices $x, y \in D$.

\begin{claim}
\label{c:claim3}
$d_G(x,y) \le 2$ for every two vertices $x, y \in D$.
\end{claim}
\proof Suppose, to the contrary, that the $\gamma$-set $D$ of $G$ contains two vertices at distance~$3$ in $G$. Let $v_1$ and $v_2$ be two vertices that belong to $D$ at distance~$3$ in $G$. 

Suppose that both $v_1$ and $v_2$ have two nonadjacent neighbors. Let $x_i$ and $y_i$ be two neighbors of $v_i$ that are not adjacent for $i \in [2]$. Let $X_1 = \{x_1,y_1\}$ and $X_2 = \{x_2,y_2\}$. Since $d_G(v_1,v_2) = 3$, we note that $X_1 \cap X_2 = \emptyset$. We now consider the induced subgraph $H = G[\{v_1,x_1,y_1,v_2,x_2,y_2\}]$. If there is no edge between $X_1$ and $X_2$, then $H = H_7$. If there is exactly one edge between $X_1$ and $X_2$, then $H = H_8$. If there are exactly two edges between $X_1$ and $X_2$, then $H = H_9$ or $H = H_{10}$. If there are exactly three edges between $X_1$ and $X_2$, then $H = H_5$. If there are four edges between $X_1$ and $X_2$, then $H = H_6$. In all cases, we contradict our assumption that $G$ is $\cH$-free.

Hence, the closed neighborhood of $v_1$ or $v_2$ (or both) induces a complete subgraph of $G$. By symmetry, and renaming vertices if necessary, we may assume that $G[N_G[v_1]]$ is a complete subgraphs of $G$. Let $v_1 w_1 w_2 v_2$ be a $(v_1,v_2)$-path in $G$. Replacing the vertex $v_1$ in the set $D$ by the vertex $w_1$ produces a new $\gamma$-set $D' = (D \setminus \{v_1\}) \cup \{w_1\}$ of $G$ noting that $N_G[v_1] \subset N_G[w_1]$ and $|D'| = |D| = \gamma(G)$. However, since the two neighbors $v_1$ and $w_2$ of $w_1$ are not adjacent, the closed neighborhood of $w_1$ does not induce a complete subgraph of $G$. Therefore the $\gamma$-set $D'$ contains fewer vertices whose closed neighborhood induces a complete subgraph of $G$ than does the set $D$, contradicting our choice of the set $D$.~\smallqed

\medskip
By Claim~\ref{c:claim3}, $d_G(x,y) \le 2$ for every two vertices $x, y \in D$. For notational simplicity, in what follows in the proof of the Theorem~\ref{gamma-alpha_c-perfect-graphs}, we let $V = V(G)$ and we simply write $N[v]$ rather than $N_G[v]$ for a vertex $v \in V$. Further for notational simplicity, if $x$ and $y$ are two vertices, then we let $G_{x,y} = G[N[x] \setminus N[y]]$.

\begin{claim}
\label{c:claim4}
There are vertices $u$ and $v$ in $D$ such that neither $G_{u,v}$ nor $G_{v,u}$ is a complete graph.
\end{claim}
\proof Suppose, to the contrary, that for every two vertices $u$ and $v$ belonging to $D$ at least one of the subgraphs $G_{u,v}$ or $G_{v,u}$ is a complete graph. Renaming vertices if necessary, we may assume that $v \in D$ and $G_{u,v}$ is a complete graph for every vertex $u \in D \setminus \{v\}$. Let $I_v$ be a maximum independent set of $G$ that contains the vertex~$v$. We note that $|I_v \cap N_G[v]| = |\{v\}| = 1$. Since $D$ is a dominating set of $G$, we have
\[
V = \bigcup_{u \in D} N[u].
\]
Therefore,
\[
\begin{array}{rcl}
I_v = I_v \cap V &=& I_v\cap \big(N[v]\cup (V\setminus N[v])\big) \1 \\
&=& (I_v\cap N[v])\cup (I_v\cap (V\setminus N[v])) \1 \\
&=& \displaystyle{ \{v\} \cup \bigg(I_v\cap \bigg(\bigcup_{u\in D\setminus\{v\}}N[u]\setminus N[v]\bigg)\bigg) } \1 \\
&=& \displaystyle{ \{v\} \cup \bigg( \bigcup_{u\in D\setminus\{v\}} (I_v \cap (N[u]\setminus N[v])) \bigg). }
\end{array}
\]

Since $I_v$ is an independent set and $G_{u,v} = G[N[u]\setminus N[v]]$ is a complete graph, we note that $|I_v \cap (N[u] \setminus N[v])| \le 1$. Hence by our earlier observations, we have
\[
\begin{array}{rcl}
\alpha_c(G) \le |I_v| &=& \displaystyle{ 1+ \left| \bigcup_{u\in D\setminus\{v\}} (I_v \cap (N[u]\setminus N[v])) \right| } \2 \\
&\le & \displaystyle{ 1+ \sum_{u\in D\setminus\{v\}} |I_v\cap (N[u]\setminus N[v])| } \2 \\ &\le & 1+ |D\setminus\{v\}| \1 \\
& = & |D| = \gamma(G) < \alpha_c(G),
\end{array}
\]
a contradiction.~\smallqed

\medskip
By Claim~\ref{c:claim4}, there are vertices $u$ and $v$ in $D$ such that neither $G_{u,v}$ nor $G_{v,u}$ is a complete graph. Let $v_1$ and $v_2$ be vertices in $N[v] \setminus N[u]$ that are not adjacent, and let $u_1$ and $u_2$ be vertices in $N[u] \setminus N[v]$ that are not adjacent. Let $A = \{v_1,v_2\}$ and $B = \{u_1,u_2\}$. We now consider the induced subgraph $H = G[\{v,v_1,v_2,u,u_1,u_2\}]$.

Suppose that $u$ and $v$ are not adjacent. If there is no edge between $A$ and $B$, then $H = H_7$. If there is exactly one edge between $A$ and $B$, then $H = H_8$. If there are exactly two edges between $A$ and $B$, then $H = H_9$ or $H = H_{10}$. If there are exactly three edges between $A$ and $B$, then $H = H_5$. If there are four edges between $A$ and $B$, then $H = H_6$. In all cases, we contradict our assumption that $G$ is $\cH$-free.

Hence, $u$ and $v$ are adjacent. Now, depending on the possible edges joining the vertices between $A$ and $B$, the induced subgraph $H$ is one of the forbidden graphs $H_1, \ldots, H_6$ shown in Figure~\ref{fig:forbidden-graphs}, once again contradicting our assumption that $G$ is $\cH$-free. This completes the proof of Theorem~\ref{gamma-alpha_c-perfect-graphs}.~\QED

\section{Consequence of Theorem~\ref{gamma-alpha_c-perfect-graphs}}

In this section, we present some consequences of our main result, namely Theorem~\ref{gamma-alpha_c-perfect-graphs}.

In 1979 Sumner and Moore~\cite{SuMo-79} have observed that a graph $G$ is $(\gamma,i)$-perfect if and only if $\gamma(H) = i(H)$ for every induced subgraph $H$ of $G$ with $\gamma(H)=2$. As a consequence of Theorem~\ref{gamma-alpha_c-perfect-graphs} a similar property holds for $(\gamma,\alpha_c)$-perfect graphs noting by Observation~\ref{ob:family} that every graph $H$ in the forbidden family $\cH$ satisfies $\gamma(H) = 2$.

\begin{corollary}
A graph $G$ is $(\gamma,\alpha_c)$-perfect if and only if $\gamma(H)= \alpha_c(H)$ for every induced subgraph $H$ of $G$ with $\gamma(H)=2$. Equivalently, a graph $G$ is $(\gamma,\alpha_c)$-perfect if and only if no induced subgraph $H$ of $G$ with $\gamma(H)=2$ has $\alpha_c(H)=3$.
\end{corollary}

\begin{figure}[htb]
\begin{center}
\begin{tikzpicture}[scale=.8,style=thick,x=0.8cm,y=0.8cm]
\def\vr{2.75pt} 
\path (0,0.5) coordinate (a1);
\path (1,0.5) coordinate (a2);
\path (1.5,1.5) coordinate (a);
\path (2,0.5) coordinate (a3);
\path (3,0.5) coordinate (a4);
\draw (a1)--(a)--(a2);
\draw (a3)--(a)--(a4);
\draw (a) [fill=white] circle (\vr);
\draw (a1) [fill=white] circle (\vr);
\draw (a2) [fill=white] circle (\vr);
\draw (a3) [fill=white] circle (\vr);
\draw (a4) [fill=white] circle (\vr);
\draw (1.5,-0.8) node {{\small (a) A star $K_{1,4}$}};
\path (5,0) coordinate (v3);
\path (5,1) coordinate (v4);
\path (6,0) coordinate (v5);
\path (6,1) coordinate (v6);
\path (6.5,2) coordinate (v7);
\path (7,1) coordinate (v8);
\path (8,1) coordinate (v9);
\path (7,0) coordinate (v1);
\path (8,0) coordinate (v2);
\draw (v7)--(v4)--(v3);
\draw (v5)--(v6)--(v7)--(v8);
\draw (v9)--(v7);
\draw (v1)--(v8);
\draw (v2)--(v9);
\draw (v1) [fill=white] circle (\vr);
\draw (v2) [fill=white] circle (\vr);
\draw (v3) [fill=white] circle (\vr);
\draw (v4) [fill=white] circle (\vr);
\draw (v5) [fill=white] circle (\vr);
\draw (v6) [fill=white] circle (\vr);
\draw (v7) [fill=white] circle (\vr);
\draw (v8) [fill=white] circle (\vr);
\draw (v9) [fill=white] circle (\vr);
\draw (6.5,-0.8) node {{\small (b) A spider}};
\path (10,0) coordinate (u3);
\path (10,1) coordinate (u4);
\path (11,0) coordinate (u5);
\path (11,1) coordinate (u6);
\path (11.5,2) coordinate (u7);
\path (12,1) coordinate (u8);
\path (13,1) coordinate (u9);
%
\draw (u7)--(u4)--(u3);
\draw (u5)--(u6)--(u7)--(u8);
\draw (u9)--(u7);
\draw (u3) [fill=white] circle (\vr);
\draw (u4) [fill=white] circle (\vr);
\draw (u5) [fill=white] circle (\vr);
\draw (u6) [fill=white] circle (\vr);
\draw (u7) [fill=white] circle (\vr);
\draw (u8) [fill=white] circle (\vr);
\draw (u9) [fill=white] circle (\vr);
\draw (11.5,-0.8) node {{\small (c) A  wounded spider}};
\path (15,1) coordinate (u3);
\path (16,1) coordinate (u4);
\path (17,1) coordinate (u5);
\path (18,1) coordinate (u6);
\path (19,0.25) coordinate (u7);
\path (19,0.75) coordinate (u8);
\path (19,1.25) coordinate (u9);
\path (19,1.75) coordinate (u10);
%
\draw (u3)--(u4)--(u5)--(u6);
\draw (u10)--(u6)--(u9);
\draw (u8)--(u6)--(u7);
\draw (u3) [fill=white] circle (\vr);
\draw (u4) [fill=white] circle (\vr);
\draw (u5) [fill=white] circle (\vr);
\draw (u6) [fill=white] circle (\vr);
\draw (u7) [fill=white] circle (\vr);
\draw (u8) [fill=white] circle (\vr);
\draw (u9) [fill=white] circle (\vr);
\draw (u10) [fill=white] circle (\vr);
\draw (16.75,-0.8) node {{\small (d) A  broom}};
\end{tikzpicture}
\caption{$(\gamma,\alpha_c)$-perfect trees}
\label{woundedspider}
\end{center}
\end{figure}
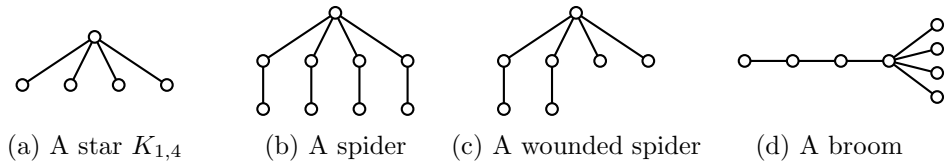

A \emph{nontrivial star} is a star $K_{1,k}$, where $k \ge 1$. A \emph{spider} is the tree obtained from a~nontrivial star by subdividing every edge exactly once, and a~\emph{wounded spider} is a tree obtained from a~nontrivial star $K_{1,k}$ by subdividing at most $k-1$ of its edges exactly once. In particular, a~nontrivial star is a wounded spider. A \emph{broom} with a handle of length~$3$ is the tree obtained from a nontrivial star $K_{1,k}$, where $k \ge 2$, by subdividing one edge twice. Examples of a~star, spider, wounded spider, and broom are illustrated in Figure~\ref{woundedspider}.

As a further consequence of Theorem \ref{gamma-alpha_c-perfect-graphs} we have the following characterization of $(\gamma,\alpha_c)$-perfect trees.

\begin{corollary}
If $G$ is a tree, then the following statements are equivalent: \\ [-24pt]
\begin{enumerate}
\item[$(1)$]  $G$ is $(\gamma,\alpha_c)$-perfect.
\item[$(2)$] $G$ does not contain any of the graphs $H_1$, $H_7$, and $H_8$ as an induced subgraph.
\item[$(3)$] $G$ has diameter at most~$4$ and $G$ has at most one vertex of degree at least~$3$.
\item[$(4)$] $G$ is $K_1$, a spider, a wounded spider, or a broom with a handle of length~$3$.
\end{enumerate}
\end{corollary}

Sumner~\cite{Su-90} showed that a chordal graph $G$ is $(\gamma,i)$-perfect if and only if $G$ is $H_1$-free. Now we have the following characterization of $(\gamma,\alpha_c)$-perfect chordal graphs and, in particular, a~characterization of $(\gamma,\alpha_c)$-perfect block graphs.

\begin{corollary}
A chordal graph $G$ is a $(\gamma,\alpha_c)$-perfect graph if and only if $G$ does not contain any of the graphs $H_1$, $H_7$, and $H_8$ as an induced subgraph. In particular, a block graph $G$ is a $(\gamma,\alpha_c)$-perfect graph if and only if it has one of the following properties:  \\ [-22pt]
\begin{enumerate}
\item[$(1)$] $G$ is a block graph of diameter~$0$,~$1$, or~$2$;
\item [$(2)$] $G$ is a block graph of diameter~$3$ with at most one vertex belonging to at least three blocks;
\item [$(3)$] $G$ is a block graph of diameter~$4$ and $G$ has one of the following properties:
\begin{enumerate}
\item[{\rm (a)}] Each inner block of $G$ is of order~$2$ and $G$ has at most one vertex belonging to at least three blocks;
\item[{\rm (b)}] $G$ has exactly one inner block of order at least~$3$ and either no vertex belonging to at least three blocks or exactly one vertex belonging to at least three blocks and this vertex belongs to the unique inner block of order at least~$3$ in $G$.
\end{enumerate}
\end{enumerate}
\end{corollary}

Recall as remarked earlier, every claw-free graph (and, therefore, every line graph) is a~$(\gamma,i)$-perfect graph.  Each of the graphs $H_7$, $H_8$, and $H_{9}$ shows that not every claw-free graph is a $(\gamma,\alpha_c)$-perfect graph. However, from Theorem~\ref{gamma-alpha_c-perfect-graphs} we have the following characterization of claw-free graphs that are $(\gamma,\alpha_c)$-perfect.

\begin{corollary}
\label{claw-free}
A claw-free graph $G$ is a $(\gamma,\alpha_c)$-perfect graph if and only if none of the graphs $H_7$, $H_8$, and $H_{9}$ is an induced subgraph of $G$.
\end{corollary}
\proof If $G$ is a $(\gamma,\alpha_c)$-perfect graph, then by Theorem~\ref{gamma-alpha_c-perfect-graphs} none of the graphs $H_7$, $H_8$, and $H_{9}$ is an induced subgraph of $G$. On the other hand, suppose that $G$ is a claw-free graph and none of the graphs $H_7$, $H_8$, and $H_{9}$ is an induced subgraph of $G$. Since every graph in the family $\cH \setminus \{H_7,H_8,H_9\}$ contains a claw, this implies that $G$ is $\cH$-free.  Consequently, $G$ is a $(\gamma,\alpha_c)$-perfect graph by Theorem~\ref{gamma-alpha_c-perfect-graphs}.~\QED

\medskip
As a consequence of Corollary~\ref{claw-free} we have the following result.

\begin{corollary}
\label{line-graphs}
A line graph $G$ is a $(\gamma,\alpha_c)$-perfect graph if and only if none of the graphs $H_7=2P_3$, $H_8=P_6$, and $H_9=C_6$ is an induced subgraph of $G$. Equivalently, if $G=L(H)$ for some graph $H$, then $G$ is a $(\gamma,\alpha_c)$-perfect graph if and only none of the graphs $2P_4$, $P_7$, and $C_6$ is a subgraph of~$H$.
\end{corollary}
\proof Since every line graph is claw-free, the first statement is immediate from Corollary~\ref{claw-free}. The equivalence is also obvious, since $2P_3$ ($P_6$, $C_6$, resp.) is an induced subgraph of the line graphs $L(H)$ if and only $2P_4$ ($P_7$, $C_6$, resp.) is a subgraph (but not necessary induced subgraph) of $H$.~\QED

\medskip
The following corollary is immediate from Corollary~\ref{line-graphs}.

\begin{corollary}
The middle graph $M(H)=L(H\circ K_1)$ of a graph $H$ is a $(\gamma,\alpha_c)$-perfect graph if and only if $H$ has no two non-adjacent edges, that is, if and only if $H$ has at most one nontrivial component, and this nontrivial component (if any) is a star.
\end{corollary}

\section{Open questions}

In this paper we provided a characterization of common domination perfect graphs in terms of ten forbidden induced subgraphs. We close this paper with the following four open questions that we have yet to settle.\\ [-22pt]
\begin{enumerate}
\item[(1)] Which  total graphs are $(\gamma,\alpha_c)$-perfect?
\item[(2)] Which $k$-trees are $(\gamma,\alpha_c)$-perfect?
\item[(3)] Which products of graphs are $(\gamma,\alpha_c)$-perfect?
\item[(4)] Which powers of graphs are $(\gamma,\alpha_c)$-perfect?
\end{enumerate}


\end{document}